\documentclass[draft,11pt]{amsart}

\newcommand{\Ad}{\mathrm{Ad}}

\newcommand{\ad}{\mathrm{ad}}

\newcommand{\rd}{\mathrm{d}}
\newcommand{\bs}{\mathbf{s}}
\newcommand{\C}{\mathbb{C}}

\newcommand{\fcirc}{\mbox{\footnotesize$\circ$}}
\newcommand{\cB}{\mathcal{B}}

\newcommand{\cU}{\mathcal{U}}

\newcommand{\diag}{\mathrm{diag}}

\newcommand{\fg}{\mathfrak{g}}
\newcommand{\fgl}{\mathfrak{gl}}

\newcommand{\fh}{\mathfrak{h}}

\newcommand{\fk}{\mathfrak{k}}
\newcommand{\fl}{\mathfrak{l}}
\newcommand{\bsl}{\boldsymbol{\lambda}}
\newcommand{\fm}{\mathfrak{m}}

\newcommand{\tfk}{\tilde{\mathfrak{k}}}

\newcommand{\fp}{\mathfrak{p}}

\newcommand{\fsl}{\mathfrak{sl}}
\newcommand{\fso}{\mathfrak{so}}

\newcommand{\fsu}{\mathfrak{su}}

\newcommand{\GL}{\mathrm{GL}}

\newcommand{\R}{\mathbb{R}}

\newcommand{\bsk}{\boldsymbol{\kappa}}

\newcommand{\SL}{\mathrm{SL}}
\newcommand{\SO}{\mathrm{SO}}
\newcommand{\SU}{\mathrm{SU}}

\newcommand{\Lstar}{\mbox{\Large $\star$}}

\makeatletter
\newcommand{\oset}[3][-0.1ex]{%
  \mathrel{\mathop{#3}\limits^{
    \vbox to#1{\kern-2\ex@
    \hbox{$\scriptstyle#2$}\vss}}}}
\makeatother

 \newtheorem{thm}{Theorem}[section]

 \newtheorem{prop}[thm]{Proposition}
 \theoremstyle{definition}
 
\theoremstyle{rem}
 \newtheorem{rem}[thm]{Remark}
 
\numberwithin{equation}{section}

\begin{document}

%
%
%
%
%

\title[Einstein-Yang-Mills Lorentzian symmetric spaces]{Lorentzian symmetric spaces which \\
are Einstein-Yang-Mills with respect \\
to  invariant metric connections}

\author[M.\,Castrill\'on L\'opez]{M.\,Castrill\'on L\'opez}
\address{Departamento de \'Algebra, Geometr\'i a y Topolog\'i a,
Facultad de Matematicas, UCM \\
Plaza de Ciencias 3, 28040-Madrid, Spain}
\email{mcastri@mat.ucm.es}

\thanks{MCL and ERM were partially funded by MCIU (Spain)
under project no. PGC2018-098321-B-I00.}

\author[P.\,M. Gadea]{P.\,M. Gadea}
\address{Instituto de F\'\i sica Fundamental, CSIC \\ Serrano 113 bis, 28006-Madrid, Spain}
\email{p.m.gadea@csic.es}

\author[M.\,E. Rosado Mar\'\i a]{M.\,E. Rosado Mar\'\i a}
\address{Departamento de Matematica Aplicada,
Escuela T\'ecnica Superior de Arquitectura \\
Universidad Polit\'ecnica de Madrid \\
Avda.\ Juan de Herrera 4, 28040-Madrid, Spain}
\email{eugenia.rosado@upm.es}


\begin{abstract}
We classify four-di\-mensional connected simply-connected
indecomposable Lorentzian symmetric spaces $M$ with connected nontrivial iso\-tropy group furnishing solutions of the Einstein-Yang-Mills  equations. Those
solutions  with respect to some invariant metric connection $\Lambda$ in the bundle of orthonormal frames of $M$ and some diagonal metric on the holonomy algebra corresponding to $\Lambda$.
\end{abstract}

\maketitle

\section{Introduction}
\label{secone}

Levichev gave in \cite{Lev1} several instances of Lorentzian symmetric Lie groups providing
solutions of the Einstein-Yang-Mills equations. Specifically: the Minkowski spacetime
$\R^{1,3}$; $\SU(2) \times \R$, which corresponds to an Einstein static universe;
the Heraclitus spacetime $\widetilde{\SU(1,1)} \times \R$ (see \cite{Lev2}
and recall that $\SU(1,1)$ is isomorphic to $\SL(2,\R)$);
and the oscillator group (Streater \cite{Str},
Calvaruso and Van der Veken \cite{CV}, Nappi and Witten \cite{NW}).
The latter space corresponds to an isotropic electromagnetic field, and was studied (with its usual Lorentzian metrics) in the vein of the present article in \cite{DGO}.

In turn, Komrakov \cite{Kom2} gave the classification of real pseudo-Riemann\-ian pairs $(\fk,\fh)$
of codimension four. He also classified
\cite[Theorem  2]{Kom3} those pairs providing solutions of the Einstein-Maxwell equations.

A natural question then arises: Which are the four-dimensional
homogeneous Lorentzian manifolds furnishing solutions of the Einstein-Yang-Mills equations?
Since there are $76$ Komrakov's reductive Lorentzian
pairs, we restrict ourselves in the present paper to a simpler question: Which are the four-dimensional
Lorentz\-ian symmetric spaces furnishing solutions of the Einstein-Yang-Mills equations? Actually, the relevant role played by symmetric spaces as models in General Relativity provides additional motivation to this question. Thanks to Komrakov's classification we have been able to obtain the classification of such spaces of a particular kind: Those providing solutions
with respect to some $K$-invariant connection $\Lambda$ in the bundle of orthonormal frames
$SO(M)$ of $M = K/H$ and a diagonal metric on the holonomy algebra $\fl$ corresponding to $\Lambda$.
To avoid repetitions of the last phrases we shall simply write
``Einstein-Yang-Mills equations of type $(\Lambda,g^{\fl}_{\mathrm{diag}})$.''

The aim of the present paper is to find the four-dimensional connected simply-connected
indecomposable Lorentzian symmetric spaces $(M=K/H$, $g)$, with connected nontrivial isotropy group $H$,
furnishing solutions of the Ein\-stein-Yang-Mills equations of type $(\Lambda,g^{\fl}_{\mathrm{diag}})$.

As it is well-known, simply-connected indecomposable
Lorentzian symmetric spaces were classified by Cahen
and Wallach \cite{CahWal}.
We use here Komrakov's classification of pairs, which includes the Lo\-rentzian symmetric pairs, and gives more
detailed expressions of some facts of the spaces under study.

Let $\fk = \fh \oplus \fm$ be a reductive decomposition. Then, the $\Ad(H)$-invariance of the scalar product
on $\fm$ furnishes in each case the relevant metrics on $M$. The construction below implies reductive
decompositions $\tilde \fk = \fl \oplus \fm$ related to the curvature operators of
the invariant connections on $(\tfk, \fl)$, and a new metric on each Lie algebra $\tilde\fk$. We restrict
ourselves, only for the sake of brevity, to a diagonal metric on each~$\fl$.

As for the contents,
after some preliminaries in Section \ref{sectwo}, we list in Section \ref{secthr},
among the $187$ Komrakov pairs, the $76$ ones which are reductive and Lorentzian and then the $38$ ones being
moreover symmetric (see Proposition \ref{threeone}, Tables \ref{aaa} and \ref{bbb}).
In Section \ref{secfou}, we exhibit the first Einstein-Yang-Mills equation for each pair in Table~\ref{bbb}.
Then we list the $10$ pairs providing some solution of the 1st EYM eq.\
of type $(\Lambda,g^{\fl}_{\mathrm{diag}})$, and the respective values of the cosmological
and gravitational constants. The corresponding currents are null (Proposition \ref{pr32}).
In Section \ref{secfiv}, we list (Theorem \ref{th41}) the ten four-dimensional connected simply-connected
indecomposable Lorentzian symmetric spaces with nontrivial connected
isotropy group, furnishing solutions of the Einstein-Yang-Mills equations
of type $(\Lambda,g^{\fl}_{\mathrm{diag}})$.

\section{Preliminaries}
\label{sectwo}

Let $K$ be a connected Lie group, whose identity element will be denoted by $e$ and let $H$ be a closed subgroup, such that
$M=K/H$ is a reductive homogeneous manifold, with reductive decomposition
$\fk = \fh \oplus \fm$, $[\fh,\fm] \subset \fm$.
Let $P\to M$ be a $G$-structure on $M$ invariant with respect to the action of $K$ on the frame bundle of $M$. We denote by $\rho \colon H\to G$ the linear isotropy representation at the origin $o \equiv eH$.
According to \cite[vol.\ II, pp.\ 191,192]{KobNom}, given a $K$-invariant $G$-structure $P$ on $M$, there is a one-to-one correspondence between the set of $K$-invariant connections $\Lambda$ in $P$ and the set of linear mappings $\Lambda_\fm \colon \fm \to \fg$ such that
\begin{equation}\label{diagso13}
\Lambda_\fm \big(\ad(Z)(X)\big) = \ad \big(\rho_*(Z)\big)\big(\Lambda_\fm(X)\big), \qquad X \in \fm, \quad Z \in \fh,
\end{equation}
where $\rho_*$ is the Lie algebra homomorphism $\fh \to \fk$ induced by $\rho$, $\ad$ denotes the adjoint representation of $\fh$ in $\fk$ and
$\ad \,\fcirc\, \rho_*$ denotes the adjoint representation of $\fg$ in $\fg$.
The correspondence is given by $\Lambda(X) = \rho_*(X)$ if $X \in \fh$ and $\Lambda(X) = \Lambda_\fm(X)$ if $X \in \fm$.
The curvature tensor $R$ of the $K$-invariant connection
corresponding to $\Lambda_\fm$ can be expressed at $o$ by
\begin{equation}
R(X,Y)_o = [\Lambda_\fm(X), \Lambda_\fm(Y)] - \Lambda_\fm([X,Y]_\fm) - \rho_*([X,Y]_\fh),
\quad X,Y \in \fm. \label{rxyo}
\end{equation}

A reductive homogeneous pseudo-Riemannian manifold is a reductive homogeneous manifold $(M=K/H,g)$,
where $g$ is a $K$-invariant pseudo-Riemannian metric on $M$. The latter metrics are known to be in bijective
correspondence with the $\Ad(H)$-invariant nondegenerate symmetric  bilinear forms $\langle\cdot,\cdot\rangle$ on $\fm$
(see, e.g., \cite[vol.\ II, p.\ 200]{KobNom}).

We are interested in four-dimensional Lorentzian symmetric manifolds $(M=K/H,g)$ and $K$-invariant
connections in the principal bundle of orthonormal frames $\mathrm{SO}(1,3)$
$\hookrightarrow SO(M) \to M$.

We now recall that Yang-Mills fields are (see, e.g., Binz, \'Sniatycki and Fischer \cite[p.\ 365]{BinSniFis})
connections in a principal bundle over a spacetime manifold $(M,g)$ satisfying a set of partial differential equations given below.
In order to give as few definitions and general results as possible, we refer to
those in Bleecker \cite[\S\S9.3]{Ble}, mainly Definitions $9.3.1$, $9.3.4$ and Theorem $9.3.3$.
We shall specialize them to the case of invariant
connections in the bundle of orthonormal frames.

In this context, the Einstein-Yang-Mills equations are the variational equations defined by the Lagrangian
density on metrics $g$ and linear connections $\Lambda$
\[
L(g,\Lambda)=\mathbf{s}v_g+R \,\dot{\wedge}\, \Lstar R,
\]
where $\mathbf{s}$ denotes the scalar curvature, $v_g$ is the volume form and $\Lstar$ the Hodge star operator defined by $g$, and $\dot{\wedge}$ is the wedge product of $\mathrm{End}(\mathfrak{m})$-valued forms defined by a convenient metric on $\mathfrak{m}$ as it follows. We consider the holonomy algebra $\fl\subset \mathfrak{o}(1,3)$ of $\Lambda$ as the Lie subalgebra of skew-symmetric
endomorphisms of $(\fm,g_o)$ generated by the operators $R_{XY}$, where $X,Y \in \fm$. We take on $\tilde{\mathfrak{k}}=\mathfrak{l}\oplus \mathfrak{m}$ an $\mathrm{Ad}(\widetilde{K})$-invariant metric $\tilde{g}$, $\widetilde{K}$ being the connected and simply connected Lie group with Lie algebra $\tilde{\fk}$. As explained above, we shall chose $\tilde{g}$ in each case being of the type
\begin{equation}
\tilde g \equiv
\begin{pmatrix}
    g_{ij}     & \mathbf{0}    \\ \noalign{\medskip}
    \mathbf{0} & g_{\alpha\alpha} \\
\end{pmatrix}, \qquad \alpha = 5,\dotsc,4 + \dim\fl,
\end{equation}
where $(g_{ij})$ stands for the $\Ad(H)$-invariant Lorentzian metric on $M$,
and $g_{\alpha,\alpha}$ is a diagonal metric for the basis of $\mathfrak{m}$,
obtained in each case from Komrakov's data.

The equations derived from this variational problem split into two. The {\it first
Einstein-Yang-Mills equation\/} is the equation defined by $g$, which reads
\begin{equation}
\label{eym1}
\mathbf{r} + \Big(\bsl  - \frac12 \bs\Big) g  = \bsk   T,
\end{equation}
where $\mathbf{r}$ and $\bs$ are the Ricci tensor and the scalar curvature of the Levi-Civita connection of $(M,g)$,
and $\bsl $, $\bsk$ stand for the cosmological and gravitational  constants. The tensor $T$ denotes the energy-momentum tensor with respect to $\Lambda$ and its local expression is given in curly brackets in \eqref{eym11} below.

\begin{rem} \em
\label{ntemt}
We shall consider as solutions of the 1st EYM eq.\  only those having nontrivial energy-momentum tensor $T$, since for $T\equiv0$, the 1st EYM eq.\ is equivalent to the manifold being Einstein.
\end{rem}

Note that, whereas the  right-hand side of the equation \eqref{eym11} depends only on $g$, the right-hand side depends on $g$, $\tilde{g}$ and $\Lambda$. The tensor components of the first EYM equation (cf.\ Bleecker \cite[Theorem  9.3.3]{Ble}, Levichev \cite[(1),(2)]{Lev1}) with respect to the basis $\{u_1,\ldots, u_4\}$ of $\mathfrak{m}$ used in Komrakov's list are
\begin{align}
  &\mathbf{r}_{ij} +\Big( \bsl  - \frac12 \sum_{k,l=1}^4  g^{kl}\mathbf{r}_{kl}\Big) g_{ij}  \label{eym11}\\
  & \qquad= \bsk  \Bigg\{ \dfrac12 \sum^4_{k,l,h,m=1}\sum^{4+\dim\fl}_{\alpha,\beta=5}
\Big( R^\alpha_{ik}R^{\beta}_{jl}g^{kl} g_{\alpha\beta} - \frac14 g_{ij}
R^{\alpha}_{kl}R^{\beta}_{hm}g^{kh}g^{lm} g_{\alpha\beta}\Big) \Bigg\}, \notag
\end{align}
for each $i\leq j$, with $i,j=1,\dotsc,4$; some $\bsl ,\bsk  \in \R$, and where $g_{ij} = g(u_i,u_j)$, $\mathbf{r}_{ij} = \mathbf{r}(u_i,u_j)$, and so on. As explained above, we shall  actually take $\beta=\alpha$.
Note that although the factor $\frac12$ could be absorbed
into the constant $\bsk$, Bleecker's formula preserves the coefficient $-\frac12$ of the self-action, which is
consistent with the classical one for the electromagnetic field (cf.\ e.g.\ \cite[\S\S5.2, 1st par., Def.\  5.2.1]{Ble}).

The {\it second Einstein-Yang-Mills equation\/} is the equation defined by $\Lambda$ and it reads
\begin{equation}
\label{eym2}
\Lstar D \Lstar R = 0,
\end{equation}
which is just the standard Yang-Mills equation for $\Lambda$. The Lagrangian of this variational problem can be modified to include currents $J$ for the Yang-Mills field $\Lambda$. In that case, the second equation would contain that current. Interestingly, as we will see, the latter situation will not make sense in the framework presented in this paper.

\section{Lorentzian symmetric pairs}
\label{secthr}

In order to classify the real reductive pairs $(\fk,\fh)$, Komrakov first proved \cite[p.\ 3, Lemma]{Kom1}
that if a homogeneous manifold $M = K/H$ admits an invariant pseudo-Riemannian metric
then the (Lie algebra homomorphism induced by the)
linear isotropy representation $\rho_*\colon \fh \to \fgl(\fm)$ of the pair $(\fk,\fh)$
is faithful. Moreover, there exists a basis of $\fm$
such that $\rho_*(\fh)$ lies in one of the Lie algebras $\fso(4)$, $\fso(1,3)$ and $\fso(2,2)$. He then considered the complexifications of the pairs $(\fk,\fh)$ and divided the
solution into the following steps: {\bf(i)} To find (up to conjugation) all possible
forms of the subalgebra $(\rho_*(\fh))^\C = \rho_*^\C(\fh^\C)$, which
is equivalent to classify (up to conjugation with respect to $\GL(4,\C)$) the subalgebras
$\fp$ of $\fso(4,\C)$; {\bf(ii)} For each such $\fp$, to find (up to equivalence of pairs) all complex pairs $(\fk^\C, \fh^\C)$ such that
$\rho_*^\C(\fh^\C)$ is conjugate to $\fp$; {\bf(iii)} For each such pair $(\fk^\C, \fh^\C)$, to find (up to equivalence of
pairs) all its real forms $(\fk,\fh)$.
These steps give rise to corresponding numberings and Komrakov labeled accordingly the pairs  with the notation $m.n^k$:
The letter $m$ denotes the dimension of the subalgebra $\fp$ of $\fso(4,\C)$, $n$ the number of the complex pair $(\fk^\C, \fh^\C)$ and $k$ the number of the real form of $(\fk^\C, \fh^\C)$. We restrict here to the cases related to $\fso(1,3)$.

We must first find the Lorentzian symmetric pairs.
Now, Komrakov's list includes some pairs which are not reductive. Among the reductive ones,
there are some pairs which are always either Riemannian or neutral.
Instead, some of them are never Riemannian or neutral; specifically,
$1.1^1(10(t=0))$, $2.1^2(1\!\!-\!\!6)$, $2.5^2(2\!-\!7)$, $3.2^2(1)$,  $3.3^2(1\!\!-\!\!4)$, $3.5^2(1\!\!-\!\!4)$, $6.1^3(1\!\!-\!\!3)$.
Some other pairs admit both Lorentzian and neutral invariant metrics:
{\small
\begin{equation}
\label{nl}
\begin{aligned}
  & 1.1^1(1\!\!-\!\!7), \quad    1.1^2(1\!\!-\!\!10,12(t=0)), \quad   1.1^3(1), \quad   1.1^4(1), \quad   1.4^1(5\!\!-\!\!26), \\
  & \qquad \qquad 2.4^1(1\!\!-\!\!3), \quad   3.2^2(2), \quad     3.5^1(1\!\!-\!\!4), \quad   4.1^2(1).
\end{aligned}
\end{equation}}

We have

\begin{prop}
\label{threeone}
The four-codimensional Lorentzian reductive pairs $(\fk,\fh)$
are the seventy-six pairs given in \emph{Table \ref{aaa}}.
The four-codimensional Lorentzian symmetric pairs $(\fk,\fh)$
are the thirty-eigth pairs listed in \emph{Table \ref{bbb}}.
\end{prop}

\begin{table}[htb]
\caption{Komrakov's Lorentzian reductive pairs.}
\label{aaa}
\begin{tabular}{|l|l|c||l|l|c|} \hline
              &  {\bf Cases}      & {\bf Lorentzian}  & & {\bf Cases}      & {\bf Lorentzian} \\ \hline \hline
${\bf 1.1}^1$ & $1\!\!-\!\!7$; $10(t=0)$  & $bd> c^2 $        & ${\bf 2.5}^2$ & $2\!\!-\!\!7$     & $a\neq 0$  \\
${\bf 1.1}^2$ & $1\!\!-\!\!10$; $12(t=0)$ & $d^2 > bc$        & ${\bf 3.2}^2$ & $1,2$             & $a\neq 0$    \\
${\bf 1.1}^3$ & $1$               & $ab\neq0$        & ${\bf 3.3}^2$ & $1\!\!-\!\!4$     & $a \neq 0$ \\
${\bf 1.1}^4$ & $1$               & $ab\neq0$        & ${\bf 3.5}^1$ & $1\!\!-\!\!4$     & $ab > 0$    \\
${\bf 1.4}^1$ & $5\!\!-\!\!26$ & $ad > 0$  & ${\bf 3.5}^2$ & $1\!\!-\!\!4$     & $ab < 0$     \\
${\bf 2.1}^2$ & $1\!\!-\!\!6$     & $ab\neq 0$        & ${\bf 4.1}^2$ & $1$              & $a \neq 0$    \\
${\bf 2.4}^1$ & $1\!\!-\!\!3$     & $ab> 0$         & ${\bf 6.1}^3$ & $1\!\!-\!\!3$             & $a \neq 0$  \\
\hline
\end{tabular}
\end{table}

\noindent{\it Proof.}
To know the conditions for the pairs in \eqref{nl} to be Lorentzian, we  consider
a basis $\cB = \{ e_1,\dotsc, e_n$, $u_1, \dotsc , u_4\}$ of $\fk$, where
$\{ e_1,\dotsc, e_n\}$ is a basis of $\fh$ and $\cU =  \{u_1,\dotsc, u_4\}$ a basis of $\fm$. We can
identify any bilinear form on $\fm$ with its matrix with respect to the basis $\cU$.
Then (see, e.g., Komrakov \cite[p.\ 122]{Kom3}),
the matrix of the invariant metric $g$ for each Komrakov's case  $n.m^k$  must satisfy the matricial equations
\begin{equation}
\label{tadad}
{}^t(\ad_{e_i})\cdot g + g \cdot(\ad_{e_i})= 0, \qquad i=1,\dotsc,n.
\end{equation}
Note that for each of these cases,  among all invariant metrics (corresponding
to nondegenerate matrices $g$ satisfying the above conditions), the Lorentzian ones
are characterized by  the  condition $\det g<0$, as  in both the Riemannian
and neutral cases one has $\det g>0$.
Starting in each case with the spaces either just before  \eqref{nl} or in \eqref{nl} with a generic symmetric matrix
and solving the $n$ equations \eqref{tadad}, one gets the matrices below, which only
depend on the type $n.m^k$, and correspond with the Lorentzian pairs  we are
looking for, whenever the conditions after a boldface $L$ are satisfied.

\smallskip

\noindent{The cases ${\bf 1.1}^1(1\!-\!7)$; and ${\bf 1.1}^1(10)\;(t=0)$:}
\[
g \equiv {\small \begin{pmatrix}
    0&0&a&0\\
    0&b&0&c\\
    a&0&0&0\\
    0&c&0&d
    \end{pmatrix}}, \quad \det g = a^2(c^2 - bd), \quad {\bf L}\; (bd > c^2).
\]
\noindent{{The cases ${\bf 1.1}^2(1\!\!-\!\!10)$; and ${\bf 1.1}^2(12)\;(t=0)$:}
\[
g \equiv  {\small \begin{pmatrix}
    a&0&0&0\\
    0&b&0&d\\
    0&0&a&0\\
    0&d&0&c
    \end{pmatrix}}, \; \det g = a^2(bc -d^2), \;\; {\bf L}\; (d^2 > bc).
\]
\noindent{{The case ${\bf 1.1}^3(1)$:}
\quad $g \equiv  {\small \begin{pmatrix}
    0&0&a&0\\
    0&b&0&0\\
    a&0&0&0\\
    0&0&0&b
    \end{pmatrix}}, \quad \det g = -a^2b^2, \quad {\bf L}\; (ab \neq 0).$

\noindent{{The case ${\bf 1.1}^4(1)$:}\quad $
g \equiv  {\small \begin{pmatrix}
    a&0&0&0\\
    0&0&0&b\\
    0&0&a&0\\
    0&b&0&0
    \end{pmatrix}}, \quad \det g = -a^2b^2, \quad {\bf L}\; (ab \neq 0).$

\noindent{The case ${\bf 1.4}^1(5\!\!-\!\!26)$:} \quad $
g \equiv  {\small \begin{pmatrix}
    0&0&-a&0\\
    0&a&0&0\\
    -a&0&b&c\\
    0&0&c&d
    \end{pmatrix}}, \quad \det g = -a^3d, \;\; {\bf L}\; (ad > 0).$
\noindent{The case ${\bf 2.1}^2(1\!\!-\!\!6)$:}\quad $
g \equiv  {\small \begin{pmatrix}
    0&0&a&0\\
    0&b&0&0\\
    a&0&0&0\\
    0&0&0&b
    \end{pmatrix}}, \quad \det g = -a^2b^2,
\quad {\bf L} \; (ab\neq0).$
\noindent{The case ${\bf 2.4}^1(1\!\!-\!\!3)$:}\quad $
g \equiv {\small \begin{pmatrix}
    0&0&-a&0\\
    0&a&0&0\\
    -a&0&0&0\\
    0&0&0&b
    \end{pmatrix}}, \quad  \det g = -a^3b, \quad {\bf L} \; (ab>0).$
\noindent{{The case ${\bf 2.5}^2(2\!\!-\!\!7)$:}
\quad $g \equiv
{\small \begin{pmatrix}
    0&0&a&0\\
    0&a&0&0\\
    a&0&b&0\\
    0&0&0&a
    \end{pmatrix}}, \quad \det g = -a^4, \quad {\bf L} \; (a \neq 0).$

\noindent{The case ${\bf 3.2}^2(1, 2)$:}\quad $
g \equiv  {\small \begin{pmatrix}
    0&0&a&0\\
    0&a&0&0\\
    a&0&0&0\\
    0&0&0&a
    \end{pmatrix}}, \quad \det g = -a^4, \quad {\bf L} \; (a\neq 0).$

\noindent{{The case ${\bf 3.3}^2(1\!\!-\!\!4)$:}\quad $
g \equiv  {\small \begin{pmatrix}
    0&0&a&0\\
    0&a&0&0\\
    a&0&b&0\\
    0&0&0&a
    \end{pmatrix}}, \quad \det g = -a^4, \quad {\bf L} \; (a\neq 0).$

\noindent{{The case ${\bf 3.5}^1(1\!\!-\!\!4)$:}\quad $
g \equiv  {\small \begin{pmatrix}
    0&0&2a&0\\
    0&a&0&0\\
    2a&0&0&0\\
    0&0&0&b
    \end{pmatrix}}, \quad \det g = -4a^3b, \quad {\bf L}\; (ab>0).$

\noindent{{The case ${\bf 3.5}^2(1\!\!-\!\!4)$:}
\quad $g \equiv \diag(a,a,a,b), \;\, \det g = a^3b, \quad {\bf L}\; (ab<0)$.

\noindent{{The case ${\bf 4.1}^2(1)$:}\quad $
g \equiv  {\small \begin{pmatrix}
    0&0&a&0\\
    0&a&0&0\\
    a&0&0&0\\
    0&0&0&a
    \end{pmatrix}},  \quad \det g = -a^4, \quad {\bf L}\; (a\neq 0)$.

\noindent{{The case ${\bf 6.1}^3(1,2,3)$:}\quad
$g \equiv \diag(a,a,a,-a), \quad \det g = -a^4, \quad {\bf L}\; (a\neq 0)$.

\smallskip

\noindent Table \ref{bbb} then follows from Table \ref{aaa} and the brackets in Komrakov \cite{Kom3}.
\begin{table}[htb]
\caption{Komrakov's Lorentzian symmetric pairs.}
\label{bbb}
\begin{tabular}{|l|l|c||l|l|c|} \hline
& {\bf Cases} & {\bf Lorentzian} & & {\bf Cases} & {\bf Lorentzian} \\ \hline \hline
${\bf 1.1}^1$ &  $7$; and $10\,(t=0)$ &  $bd > c^2$ & ${\bf 2.5}^2$ & $4\!\!-\!\!7$ &  $a \neq 0$ \\
${\bf 1.1}^2$ &  $9,10$; and $12\,(t=0)$ &  $d^2 > bc$ & ${\bf 3.2}^2$ & $2$ &  $a\neq 0$ \\
${\bf 1.1}^3$ & $1$ & $ab\neq 0$ & ${\bf 3.3}^2$ & $2,3,4$ & $a \neq 0$ \\
${\bf 1.1}^4$ & $1$ & $ab\neq 0$ & ${\bf 3.5}^1$ & $2,3,4$ & $ab > 0$ \\
${\bf 1.4}^1$ & $24,25,26$ & $ad > 0$ & ${\bf 3.5}^2$ & $2,3,4$ & $ab < 0$ \\
${\bf 2.1}^2$ & $1\!\!-\!\!6$  & $ab\neq 0$ & ${\bf 4.1}^2$ & $1$ & $a\neq 0$ \\
$ {\bf 2.4}^1$ & $3$ & $ab > 0$ & ${\bf 6.1}^3$ & $1,2,3$ & $a \neq 0$   \\
\hline
\end{tabular}
\end{table}
\hfill $\square$

We shall omit, for the sake of brevity, some data from Komrakov \cite{Kom3},
a paper which is easily available online.
Since the metric is invariant, the coefficients of the Levi-Civita connection
$\nabla$ can be calculated by the (short) Koszul formula
$2 g(\nabla_{u_i}u_j,u_k) = g([u_i,u_j],u_k) - g([u_j,u_k],u_i) + g([u_k,u_i],u_j)$,
and from them, the Einstein tensor can be derived.
In the next section we shall write, for the cases in Table \ref{bbb}, some facts which either are included in Komrakov's tables,
or can be calculated starting from them. Specifically:
The pair corresponding with each case (where sometimes it appears as a factor
the usual Heisenberg Lie algebra $\fh_3$ or the five-dimensional Heisenberg
Lie algebra $\fh_5$ given by the generators and nonzero brackets:
$\fh_3 = \langle h,p_1,p_2\rangle$, $[p_1,p_2] = -[p_2,p_1] = h$ and
$\fh_5 = \langle h,p_1,p_2,q_1,q_2\rangle$, $[p_1,q_1] = -[q_1,p_1] = h$, $[q_1,q_2] = -[q_2,q_1] = h)$;
the nonzero brackets; the metric and the Ricci tensor; the scalar curvature;
the matrices $\rho_*(e_i)$, $i=1,\dotsc, \dim \fl$, and $\Lambda_\fm(u_j)$, $j=1,\dotsc,4$; and the 1st EYM  eq., with its solutions, if any.

\section{The first Einstein-Yang-Mills equation}
\label{secfou}

\subsection{The cases ${\bf 1.1^1}(7)$; and ${\bf 1.1^1}(10)\;(t=0)$.}
The corresponding pairs are $\big(\fsl(2,\R) \oplus  \R^2, \; \fso(1,1)\big)$
and
{\small
\[
\left( \R^4 \times
\left\{\begin{pmatrix}
    x &0&0&0\\
    0&0 &0&0\\
    0&0&-x &0\\
    0&0&0&0
\end{pmatrix}    \right\} ,
\left\{
\begin{pmatrix}
    x &0&0&0\\
    0&0 &0&0\\
    0&0&-x &0\\
    0&0&0&0
    \end{pmatrix}
\right\}\right), \quad x \in\R.
\]}

\noindent The common
nonzero brackets are $[e_1,u_1] = u_1$, $[e_1,u_3] = -u_3$,
and the different ones are $[u_1,u_3] = e_1$ for $1.1^1(7)$ and $[u_1,u_3] = 0 $ for $1.1^1(10)\;(t=0)$.
In the first case one has
{\small
\begin{gather*}
g \equiv  {\setlength\arraycolsep{4pt}
    \begin{pmatrix}
    0&0&a&0\\
    0&b&0&c\\
    a&0&0&0\\
    0&c&0&d
    \end{pmatrix}, \;
\mathbf{r} \equiv
\begin{pmatrix}
0 & 0 & -1 & 0 \\
0 & 0&0&0 \\
-1 & 0&0&0 \\
0 & 0&0&0
\end{pmatrix}, \;\, \bs = -\dfrac2a, \;
\rho_*(e_1) =
\begin{pmatrix}
    1&0&0&0\\
    0&0&0&0\\
    0&0&-1&0\\
    0&0&0&0
    \end{pmatrix}};
\end{gather*}}

\noindent for $1.1^1(10)\;(t=0)$ one has
$\mathbf{r} \equiv {\bf 0}$, $\bs =0$, and the general relations in \eqref{diagso13} are here
$\rho_*(e_1)\cdot \Lambda_\fm(u_i) - \Lambda_\fm(u_i) \cdot \rho_*(e_1) - \Lambda_\fm([e_1,u_i]) = 0$, $i=1,\dotsc,4$.
Thus, the nonzero matrices $\Lambda_\fm(u_i)$ are
{\small \begin{gather*}
\Lambda_\fm(u_2) =
\begin{pmatrix}
    0&0&0&0\\ \noalign{\smallskip}
    0&0&0&v_{24}\\ \noalign{\smallskip}
    0&0&0&0\\ \noalign{\smallskip}
    0&-v_{24}&0&0
\end{pmatrix},  \quad \Lambda_\fm(u_4) =
\begin{pmatrix}
    0&0&0&0\\ \noalign{\smallskip}
    0&0&0&z_{24}\\ \noalign{\smallskip}
    0&0&0&0\\ \noalign{\smallskip}
    0&-z_{24}&0&0
\end{pmatrix}.
\end{gather*}}

\noindent Then, according to \eqref{rxyo}, all components of $R$ at $o$ for $1.1^1(7)$ and $1.1^1(10)$ vanish except for $R_{13} = -\rho_*(e_1)$ in the case $1.1^1(7)$. The 1st EYM equation is
{\small
\[
\begin{pmatrix}
    0&0&\bsl  a&0\\  \noalign{\smallskip}
    0&\bsl  b + (b/a)&0&\bsl  c+ (c/a)\\  \noalign{\smallskip}
    \bsl  a&0&0&0\\  \noalign{\smallskip}
    0&\bsl  c+  (c/a)&0&\bsl  d+  (d/a)
\end{pmatrix}
=
\bsk
\begin{pmatrix}
0&0&-1/2a&0\\ \noalign{\smallskip}
0&b/2a^2&0&c/2a^2\\ \noalign{\smallskip}
-1/2a&0&0&0\\ \noalign{\smallskip}
0&c/2a^2&0&d/2a^2
\end{pmatrix},
\]}

\noindent with solutions $\bsl = -1/2a$, $\bsk = a$, $a\neq 0$.

In the case $1.1^1(10)\; (t=0)$, one has $T \equiv 0$. Hence, according to Remark \ref{ntemt},
this case does not furnish any solution of the 1st EYM eq.

\subsection{The cases ${\bf 1.1^2}(9)$, ${\bf 1.1^2}(10)$; and ${\bf 1.1^2}(12)\; (t=0)$.}
The corresponding pairs are $\big(\fsu(2) \oplus  \R^2, \, \fso(2)\big)$,  $\big(\fsl(2,\R) \oplus  \R^2, \; \fso(2)\big)$ and
{\small
\[
\left( \! \R^4 \times \left\{\begin{pmatrix}
    0&0&-x&0\\
    0&0&0&0\\
    x&0&0&0\\
    0&0&0&0
\end{pmatrix}\! \right\}  ,
\left\{ \!
\begin{pmatrix}
 0&0&-x&0\\
    0&0&0&0\\
    x&0&0&0\\
    0&0&0&0
    \end{pmatrix}
\right\}\!\right),\, x\in\R,
\]}

\noindent since the metric for the case $1.1^2(12)$ is Lorentzian only for $t=0$.

The common nonzero brackets are $[e_1,u_1] = u_3$, $[e_1,u_3] = -u_1$,
and the different brackets are $[u_1,u_3] = \varepsilon e_1$, for $\varepsilon = 1,-1,0$, respectively.

In the cases $1.1^2(9,10)$ we have that
{\small
\[
g \equiv  \begin{pmatrix}
    a&0&0&0\\
    0&b&0&d\\
    0&0&a&0\\
    0&d&0&c
    \end{pmatrix}, \;
\mathbf{r} \equiv \begin{pmatrix}
\varepsilon&0&0&0\\
0&0&0&0\\
0&0&\varepsilon&0\\
0&0&0&0
\end{pmatrix}, \; \bs = \dfrac{2\varepsilon}{a},
\;
\rho_*(e_1) =
\begin{pmatrix}
0&0&-1&0\\
0&0&0&0\\
1&0&0&0\\
0&0&0&0
\end{pmatrix}.
\]}

Proceeding as above, one obtains matrices $\Lambda_\fm(u_i)$ as those
in the case $1.1^1(7)$, so that, according to \eqref{rxyo}, all components of $R$ at $o$ for
$1.1^2(9,10,12)$ vanish except for
$R_{13} =-\varepsilon \rho_*(e_1)$
in the cases $1.1^2(9)$ and $1.1^2(10)$, respectively. Then, the 1st EYM equation is
{\small
\[
\begin{pmatrix}
    \bsl  a&0&0&0\\  \noalign{\smallskip}
    0&\bsl  b -\varepsilon (b/a)&0&\bsl  d-\varepsilon (d/a)\\  \noalign{\smallskip}
    0& 0& \bsl  a&0\\  \noalign{\smallskip}
    0&\bsl  d-\varepsilon(d/a)&0&\bsl  c-\varepsilon(c/a)
\end{pmatrix}
=
\bsk
\begin{pmatrix}
    1/2a&0&0&0\\ \noalign{\smallskip}
    0&-b/2a^2&0&-d/2a^2\\ \noalign{\smallskip}
    0&0&1/2a&0\\ \noalign{\smallskip}
    0&-d/2a^2&0&-c/2a^2
\end{pmatrix},
\]}

\noindent with solutions $\bsl = \varepsilon/2a$, $\bsk = \varepsilon a$, $a \in\R\backslash\{0\}$.

According to \cite[p.\ 131]{Kom3}, the metric in the cases
$1.1^2(12(t\in [0,1]))$ is Lorentzian only for $t=0$. Then
$T \equiv 0$, so this case does not provide any solution of the 1st EYM eq.

\subsection{The case ${\bf 1.1^3}(1)$.}
The corresponding pair is
{\small
\begin{eqnarray*}
\left(\!    \R^4 \!\times\!
\left\{ \!\begin{pmatrix}
x&0&0&0 \\
0&0&0&-\lambda x\\
0&0&-x&0 \\
0&\lambda x&0&0 \\
\end{pmatrix}\!\right\},\!
\left\{\!
\begin{pmatrix}
x&0&0&0 \\
0&0&0&-\lambda x\\
0&0&-x&0 \\
0&\lambda x&0&0 \\
\end{pmatrix} \!\right\}\!\right), \; x\in \R, \, \lambda(0,1].
\end{eqnarray*}}

\noindent The nonzero brackets are $[e_1,u_1] = u_1$, $[e_1,u_2] = \lambda u_4$, $[e_1,u_3] = - u_3$, $[e_1,u_4] = - \lambda u_2$.
In this case, one has $T \equiv 0$. Hence, according to Remark \ref{ntemt},
this case does not furnish any solution of the 1st EYM eq.

\subsection{The case ${\bf 1.1^4}(1)$.}
The corresponding pair is
{\small
\begin{eqnarray*}
\left(\!    \R^4 \!\times\!
\left\{ \!\begin{pmatrix}
0&0&-x&0 \\
0&-\lambda x&0&0 \\
x&0&0&0 \\
0&0&0&\lambda x \\
\end{pmatrix}\!\right\},\!
\left\{\!
\begin{pmatrix}
0&0&-x&0 \\
0&-\lambda x&0&0 \\
x&0&0&0 \\
0&0&0&\lambda x \\
\end{pmatrix} \!\right\}\!\right), \; x\in \R, \, \lambda(0,1].
\end{eqnarray*}}

The nonzero brackets are $[e_1,u_1] = u_3$, $[e_1,u_2] = -\lambda u_2$, $[e_1,u_3] = - u_1$, $[e_1,u_4] = - \lambda u_4$.
In this case, one has $T \equiv 0$. Hence, according to Remark \ref{ntemt},
this case does not furnish any solution of the 1st EYM eq.

\subsection{The cases ${\bf 1.4^1}(24)$, ${\bf 1.4^1}(25)$ and ${\bf 1.4^1}(26)$.}
The corresponding pairs are
{\small
\begin{eqnarray*}
\left( (\fh_3 \oplus  \R) \times
\left\{
\begin{pmatrix}
0&0&0&0 \\
0&0&\varsigma x&0 \\
0&x&0&0 \\
0&0&0&0 \\
\end{pmatrix} : x \in \R \right\} , \;
\langle p \rangle\right),\quad \varsigma =+1,-1, \;\, \text{resp.,  and}\\
\left(    \R^4 \times
\left\{ \begin{pmatrix}
0&x&0&0 \\
0&0&x&0 \\
0&0&0&0 \\
0&0&0&0 \\
\end{pmatrix}\right\},
 \;
\left\{
\begin{pmatrix}
0&x&0&0 \\
0&0&x&0 \\
0&0&0&0 \\
0&0&0&0 \\
\end{pmatrix} \right\}\right), \quad x\in \R.
\end{eqnarray*}}

The common nonzero brackets are $[e_1,u_2] = u_1$, $[e_1,u_3] = u_2$,
and the different ones are $[u_2,u_3] = \varepsilon e_1$, for $\varepsilon = 1,-1,0$,
respectively.

We have in all cases $1.4^1(24,25,26)$ that
{\small \[
g \equiv  \begin{pmatrix}
    0&0&-a&0\\
    0&a&0&0\\
    -a&0&b&c\\
    0&0&c&d
    \end{pmatrix},\quad \bs =  0, \quad
\rho_*(e_1) =
   \begin{pmatrix}
    0&1&0&0\\
    0&0&1&0\\
    0&0&0&0\\
    0&0&0&0
    \end{pmatrix},
\]}

\noindent whereas $\mathbf{r} \equiv \diag(0,0,-\varepsilon,0)$, respectively.
Proceeding then as above, one obtains that $\Lambda_\fm(u_i)=0$, $i=1,\dotsc,4$,
hence the only nonzero component of $R$ at $o$ is $R_{23} = - \varepsilon \rho_*(e_1)$.

The 1st EYM eq.\ is
{\small
\[
\begin{pmatrix}
    0&0&- \bsl  a&0\\  \noalign{\smallskip}
    0&\bsl  a&0&0\\  \noalign{\smallskip}
    - \bsl  a&0&\bsl  b-\varepsilon&\bsl  c\\  \noalign{\smallskip}
    0&0        &\bsl  c&\bsl  d
\end{pmatrix}
 = \dfrac{\bsk(bd-c^2)}{4a^2d}
 {\setlength\arraycolsep{2pt}
\begin{pmatrix}
    4ad/(bd-c^2)&0&-1&0\\ \noalign{\smallskip}
    0&1&0&0\\ \noalign{\smallskip}
    -1&0&-b&-c\\ \noalign{\smallskip}
    0&0&-c&-d
\end{pmatrix},}
\]}
with no solutions for the cases $1.4^1(24,25)$.

In the case $1.4^1(26)$, one has $T \equiv 0$. Hence, according to Remark \ref{ntemt},
this case does not furnish any solution of the 1st EYM eq.

\subsection{The cases ${\bf 2.1^2}(1)$, ${\bf 2.1^2}(2)$, ${\bf 2.1^2}(3)$}
In these cases, the corresponding pairs are
$\big(   \fso(3)\oplus\fsl(2,\R), \;  \fso(2)\oplus\fso(1,1)\big)$,
$\big(\fsl(2,\R) \oplus  \fsl(2,\R), \; \fso(2)\oplus\fso(1,1)  \big)$ and
$\big((\R^2 \times \fso(2) )\oplus\fsl(2,\R), \;\fso(2)\oplus\fso(1,1)  \big)$,
respectively.
The common nonzero brackets are $[e_1,u_1] = u_1$, $[e_1,u_3]$ $= -u_3$, $[e_2,u_2]$ $= u_4$,
$[e_2,u_4] = -u_2$, $[u_1,u_3] =  e_1$, and the different ones are $[u_2,u_4] = \varepsilon e_2$,
with $\varepsilon = 1,-1,0$, respectively.

One has in all cases $2.1^2(1\!\!-\!\!6)$ that
{\small
\begin{eqnarray*}
g \equiv  \begin{pmatrix}
    0&0&a&0\\
    0&b&0&0\\
    a&0&0&0\\
    0&0&0&b
    \end{pmatrix},
\quad
\mathbf{r} \equiv
\begin{pmatrix}
0&0&-1&0 \\
0&1&0&0 \\
-1&0&0&0 \\
0&0&0&\varepsilon
\end{pmatrix}, \quad \bs = \dfrac{2(\varepsilon a-b)}{ab},
\end{eqnarray*}}
{\small
\begin{gather*}
\rho_*(e_1) = \begin{pmatrix}
    1&0&0&0\\
    0&0&0&0\\
    0&0&-1&0\\
    0&0&0&0
    \end{pmatrix}, \quad
\rho_*(e_2) = \begin{pmatrix}
    0&0&0&0\\
    0&0&0&-1\\
    0&0&0&0\\
    0&1&0&0
    \end{pmatrix}.
\end{gather*}}

The general relations in \eqref{diagso13} are here
$\rho_*(e_j)\cdot \Lambda_\fm(u_i) - \Lambda_\fm(u_i) \cdot \rho_*(e_j)
- \Lambda_\fm([e_j,u_i]) = 0$, $i=1,\dotsc,4$, $j=1,2$.
Thus in the three cases $2.1^2(1,2,3)$ the nonzero matrices $\Lambda_\fm(u_i)$ are
{\small
\begin{gather*}
\Lambda_\fm(u_2) =
\begin{pmatrix}
    0&0&0&0\\ \noalign{\smallskip}
    0&0&0&v_{24}\\ \noalign{\smallskip}
    0&0&0&0\\ \noalign{\smallskip}
    0&-v_{24}&0&0
\end{pmatrix},  \quad \Lambda_\fm(u_4) =
\begin{pmatrix}
    0&0&0&0\\ \noalign{\smallskip}
    0&0&0&z_{24}\\ \noalign{\smallskip}
    0&0&0&0\\ \noalign{\smallskip}
    0&-z_{24}&0&0
\end{pmatrix}.
\end{gather*}}

Then, the non null components of $R$ at $o$  are
$R_{13} = - \rho_*(e_1)$, $R_{24} = - \varepsilon \rho_*(e_2)$ and the 1st EYM equation is
{\small
\[ {\setlength\arraycolsep{1pt}
\begin{pmatrix}
    0&0&\bsl  a -\varepsilon(a/b)&0\\
    0&\bsl  b +(b/a)&0&0\\
    \bsl  a -\varepsilon a/b&0&0&0\\
    0&0&0& \bsl  b +(b/a)
\end{pmatrix}}
= \dfrac{\bsk((\varepsilon a)^2+b^2)}{2a^2b}
{\setlength\arraycolsep{2pt}
\begin{pmatrix}
    0&0&-1&0\\ \noalign{\smallskip}
    0&1&0&0\\ \noalign{\smallskip}
    -1&0&0&0\\ \noalign{\smallskip}
    0&0&0&1
\end{pmatrix}},
\]}

\noindent with solutions
$\bsl = (\varepsilon a-b)/2ab$,  $\bsk = ab(\varepsilon a+b)/((\varepsilon a)^2+b^2)$, $a,b\neq0$, $\varepsilon a\neq -b$.

\subsection{The cases ${\bf 2.1^2}(4)$, ${\bf 2.1^2}(5)$ and ${\bf 2.1^2}(6)$.}
The corresponding pairs are
$\big(\fso(3)  \oplus    (\R^2 \times \fso(1,1)), \fso(2)  \oplus    \fso(1,1)\big)$,
\!$\big(\fsl(2,\R)  \oplus  (\R^2 \times \fso(1,1)), \fso(2)  \oplus    \fso(1,1)\big)$ and
{\small \[
\left( \R^4 \times \left\{ \begin{pmatrix}
    x&0&0&0\\
    0&0&0&-y\\
    0&0&-x&0\\
    0&y&0&0
    \end{pmatrix}\right\} , \; \left\{ \begin{pmatrix}
    x&0&0&0\\
    0&0&0&-y\\
    0&0&-x&0\\
    0&y&0&0
    \end{pmatrix}\right\}\right), \quad x,y \in \R.
\]}

The common nonzero brackets are $[e_1,u_1] = u_1$, $[e_1,u_3] = -u_3$, $[e_2,u_2]$ $u_4$, $[e_2,u_4] = -u_2$,
and the different ones are $[u_2,u_4] = \varepsilon e_2$, for $\varepsilon = 1,-1,0$, respectively.

One has $\mathbf{r} \equiv \diag(0,\varepsilon,0, \varepsilon)$, $\bs = 2\varepsilon/b$,
{\small
\[
\rho_*(e_1) =
\begin{pmatrix}
    1&0&0&0\\ \noalign{\smallskip}
    0&0&0&0\\ \noalign{\smallskip}
    0&0&-1&0\\ \noalign{\smallskip}
    0&0&0&0
\end{pmatrix}, \quad
\rho_*(e_2) =
\begin{pmatrix}
    0&0&0&0\\ \noalign{\smallskip}
    0&0&0&-1\\ \noalign{\smallskip}
    0&0&0&0\\ \noalign{\smallskip}
    0&1&0&0
\end{pmatrix}.
\]}

The nonzero $\Lambda_\fm(u_i)$ are
{\small
\[
\Lambda_\fm(u_2)
= \begin{pmatrix}
    0&0&0&0\\ \noalign{\smallskip}
    0&0&0&v_{24}\\ \noalign{\smallskip}
    0&0&0&0\\ \noalign{\smallskip}
    0&-v_{24}&0&0
\end{pmatrix} \quad
\Lambda_\fm(u_4)
= \begin{pmatrix}
    0&0&0&0\\ \noalign{\smallskip}
    0&0&0&z_{24}\\ \noalign{\smallskip}
    0&0&0&0\\ \noalign{\smallskip}
    0&-z_{24}&0&0
\end{pmatrix}.
\]}

\noindent The only nonnull component of $R$ at $o$ is $R_{24} = - \varepsilon \rho_*(e_2)$.
The 1st EYM eq.\ is
{\small
\[
\begin{pmatrix}
    0&0&\bsl a  -\varepsilon(a/b) &0\\
    0&\bsl b&0&0\\
    \bsl  a -\varepsilon(a/b)&0&0&0\\
    0&0&0&\bsl b
\end{pmatrix}
= \dfrac{\bsk\, \varepsilon^2}{2b}
\begin{pmatrix}
    0&0&-a/b&0\\ \noalign{\smallskip}
    0&1&0&0\\ \noalign{\smallskip}
    -a/b&0&0&0\\ \noalign{\smallskip}
    0&0&0&1
\end{pmatrix}.
\]}

Hence, the cases $2.1^2(4,5)$ (but not $2.1^2(6)$) furnish the following solutions of the 1st EYM eq.: $\bsl = \varepsilon/2b$, $\bsk =  \varepsilon b$, $b\neq0$,
respectively.

\subsection{The case ${\bf 2.4^1}(3)$.}
The corresponding pair is
{\small
\[
\left( \R^4 \times \left\{
\begin{pmatrix}
x&y&0&0 \\ \noalign{\smallskip}
0&0&y&0 \\ \noalign{\smallskip}
0&0&-x&0 \\ \noalign{\smallskip}
0&0&0&0
\end{pmatrix} \right\} , \;
\left\{
\begin{pmatrix}
x&y&0&0 \\ \noalign{\smallskip}
0&0&y&0 \\ \noalign{\smallskip}
0&0&-x&0 \\ \noalign{\smallskip}
0&0&0&0
\end{pmatrix} \right\}\right), \quad x,y \in \R,
\]}

\noindent and the nonzero brackets are $[e_1,e_2] = e_2$, $[e_1,u_1] = u_1$, $[e_1,u_3] = -u_3$,
$[e_2,u_2] = u_1$, $[e_2,u_3] = u_2$.

\noindent One has $\mathbf{r} \equiv 0$, $\bs = 0$, \;
{\small
$
\rho_*(e_1) =
\begin{pmatrix}
    1&0&0&0\\
    0&0&0&0\\
    0&0&-1&0\\
    0&0&0&0
    \end{pmatrix}$, \,
$\rho_*(e_2) =
\begin{pmatrix}
    0&1&0&0\\
    0&0&1&0\\
    0&0&0&0\\
    0&0&0&0
    \end{pmatrix}.
$}

\noindent All matrices $\Lambda_\fm(u_i)$, $i=1,\dotsc,4$, are null.
Then, from the brackets of the case and \eqref{rxyo}, one has $R=0$, hence the case
${2.4^1}(3)$ does not furnish any solution of the first EYM eq.

\subsection{The cases ${\bf 2.5^2}(4)$, ${\bf 2.5^2}(5)$, ${\bf 2.5^2}(6)$ and ${\bf 2.5^2}(7)$.}
The pairs corresponding to $2.5^2(4)$ and $2.5^2(5)$ are
{\small
\[
\left( \fh_5 \times \left\{
\begin{pmatrix}
0&0&0&0&0 \\ \noalign{\smallskip}
0&0&0&\varsigma(1+t)x&0 \\ \noalign{\smallskip}
0&0&0&0&x \\ \noalign{\smallskip}
0&-x&0&0&0 \\ \noalign{\smallskip}
0&0&\varsigma(t-1)x&0&0
\end{pmatrix} : x \in \R, \; t \geq 0
\right\} , \;
\langle p_1,q_2 \rangle \right),
\]}

\noindent with $\varsigma = 1,-1$, respectively, and those for
$2.5^2(6)$ and $2.5^2(7)$ are
{\small
\begin{gather*}
\left(  \fh_5 \times \left\{
\begin{pmatrix}
0&0&0&0&0 \\ \noalign{\smallskip}
0&0&-x&0&0 \\ \noalign{\smallskip}
0&0&0&0&x \\ \noalign{\smallskip}
0&-x&0&0&0 \\ \noalign{\smallskip}
0&0&0&x&0
\end{pmatrix} : \; x \in \R
\right\} , \;\,
\fh = \langle p_1,q_2 \rangle \right), \\ \noalign{\smallskip}
\left( \R^4 \times  \left\{
\begin{pmatrix}
0 & x & 0 & -y \\
0 & 0 & -x & 0 \\
0 & 0 & 0 & 0 \\
0&0&y&0
 \end{pmatrix}\right\} , \;\left\{
\begin{pmatrix}
0 & x & 0 & -y \\
0 & 0 & -x & 0 \\
0 & 0 & 0 & 0 \\
0&0&y&0
\end{pmatrix}\right\}\right), \quad x,y\in\R.
\end{gather*}}

The common nonzero brackets in the four cases are
$[e_1,u_2] = u_1$, $[e_1,u_3] = -u_2$, $[e_2,u_3] = u_4$, $[e_2,u_4] = -u_1$,
and the nonnull remaining brackets are
$[u_2,u_3] = \varepsilon(1+t)e_1$, $[u_3,u_4] = \varepsilon(1-t)e_2$, for $\varepsilon = 1,-1,0$; and
$[u_2,u_3] = e_2$, $[u_3,u_4] = e_1$,
for $2.5^2(4)$, $2.5^2(5)$, $2.5^2(7)$; and $2.5^2(6)$, respectively.

The cases $2.5^2(4,5,6,7)$ share
{\small $g \equiv
\begin{pmatrix}
    0&0&a&0\\
    0&a&0&0\\
    a&0&b&0\\
    0&0&0&a
    \end{pmatrix}$},
$\mathbf{r} \equiv 0$, $\bs = 0$,
{\small \[
\rho_*(e_1) = \begin{pmatrix}
    0&1&0&0\\
    0&0&-1&0\\
    0&0&0&0\\
    0&0&0&0
    \end{pmatrix},\quad
\rho_*(e_2) = \begin{pmatrix}
    0&0&0&-1\\
    0&0&0&0\\
    0&0&0&0\\
    0&0&1&0
    \end{pmatrix}.
\]}

All matrices $\Lambda_\fm(u_i)$, $i=1,\dotsc,4$, vanish.
Then, according to \eqref{rxyo}, the nonnull components of $R$ at $o$ in the cases $2.5^2(4,5,7)$ are
$R_{23} =-\varepsilon (1+t)\rho_*(e_1)$,
$R_{34} = -\varepsilon (1-t)\rho_*(e_2)$
and the 1st EYM equation in the cases $2.5^2(4,5)$  is
{\small
\[
\begin{pmatrix}
0 & 0 & \bsl  a & 0 \\
0 & \bsl  a & 0 & 0 \\
\bsl  a & 0 &\bsl  b  & 0 \\
0&0&0&\bsl  a
 \end{pmatrix}
 = \dfrac{\bsk \,\varepsilon^2 b(1+t^2)}{2a^2}
\begin{pmatrix}
    4a/b&0&1&0\\ \noalign{\smallskip}
    0&1&0&0\\ \noalign{\smallskip}
    1&0&b&0\\ \noalign{\smallskip}
    0&0&0&1
\end{pmatrix}.
\]}
The nonzero components of $R$ for  $2.5^2(6)$ are
$R_{23} =-\rho_*(e_2)$, $R_{34} = -\rho_*(e_1)$.

The right-hand side in the case $2.5^2(6)$ is the same as in cases $2.5^2(4,5)$ but with $t=0$.
Hence, the cases $2.5^2(4,5,6,7)$ do not provide solutions of the 1st EYM eq.

\subsection{The case ${\bf 3.2^2}(2)$.}
The corresponding pair is
{\small
\[
\left(\! \R^4 \!\times\! \left\{\begin{pmatrix}
x & y & 0 & -z \\
0 & 0 & -y & -\lambda x \\
0 & 0 &-x & 0 \\
0&\lambda x&z&0
 \end{pmatrix}\!\right\} ,
\left\{\!\begin{pmatrix}
x & y & 0 & -z \\
0 & 0 & -y & -\lambda x \\
0 & 0 &-x & 0 \\
0&\lambda x&z&0
 \end{pmatrix}\!\right\}\!\right), \; x,y,z \in \R, \; \lambda \geq 0.
\]}

The nonzero brackets are
$[e_1,e_2] = e_2-\lambda e_3$, $[e_1,e_3] = e_3 + \lambda e_2$, $[e_1,u_1] = u_1$, $[e_1,u_2] = \lambda u_4$,
$[e_1,u_3] = -u_3$, $[e_1,u_4] = -\lambda u_2$, $[e_2,u_2] = u_1$, $[e_2,u_3] = -u_2$, $[e_3,u_3] = u_4$, $[e_3,u_4] = -u_1$.
In this case, one has $T \equiv 0$. Hence, according to Remark \ref{ntemt},
this case does not furnish any solution of the 1st EYM eq.

\subsection{The cases ${\bf 3.3^2}(2)$, ${\bf 3.3^2}(3)$ and ${\bf 3.3^2}(4)$.}
The corresponding pairs are
{\small
\[
\left( \fh_5 \times \left\{
\begin{pmatrix}
0&0&0&0&0 \\ \noalign{\smallskip}
0&0&0&\varsigma  y&x \\ \noalign{\smallskip}
0&0&0&x&y \\ \noalign{\smallskip}
0&-y&-x&y&0 \\ \noalign{\smallskip}
0&-x&-\varsigma y&0&0 \\ \noalign{\smallskip}
\end{pmatrix}
\right\} , \;
\left\langle
\begin{pmatrix}
0&0&0&0&0 \\ \noalign{\smallskip}
0&0&0&0&1 \\ \noalign{\smallskip}
0&0&0&1&0 \\ \noalign{\smallskip}
0&0&-1&0&0 \\ \noalign{\smallskip}
0&-1&0&0&0 \\ \noalign{\smallskip}
\end{pmatrix},p_1,q_2 \right\rangle \right),
\]}

\noindent for $x,y \in \R$, $\varsigma = 1,-1$, and
{\small
\[
\left( \R^4 \times \left\{\begin{pmatrix}
0 & y & 0 & -z \\
0 & 0 & -y & -x \\
0 & 0 &0 & 0 \\
0&x&z&0
 \end{pmatrix}\right\} ,\;
\left\{\begin{pmatrix}
0 & y & 0 & -z \\
0 & 0 & -y & -x \\
0 & 0 &0 & 0 \\
0&x&z&0
 \end{pmatrix}\right\}\right), \quad x,y,z \in \R,
\]}

\noindent respectively.

The common nonzero brackets are
$[e_1,e_2] = -e_3$, $[e_1,e_3] = e_2$, $[e_1,u_2] = u_4$, $[e_1,u_4] = -u_2$,
$[e_2,u_2] = u_1$, $[e_2,u_3] = -u_2$, $[e_3,u_3] = u_4$, $[e_3,u_4] = -u_1$,
and the different ones are $[u_2,u_3] = \varepsilon e_2$, $[u_3,u_4] = \varepsilon e_3$,
where $\varepsilon = 1,-1,0$, respectively.
The cases $3.3^2(2,3,4)$ share
{\small
\begin{gather*}
g \equiv  \begin{pmatrix}
    0&0&a&0\\
    0&a&0&0\\
    a&0&b&0\\
    0&0&0&a
    \end{pmatrix}, \quad
\mathbf{r} \equiv
\begin{pmatrix}
0&0&0&0 \\
0&0&0&0 \\
0&0&2\varepsilon&0 \\
0&0&0&0
\end{pmatrix}, \quad \bs = 0,
\\
\rho_*(e_1) = {\setlength\arraycolsep{3pt}
\begin{pmatrix}
    0&0&0&0\\
    0&0&0&-1\\
    0&0&0&0\\
    0&1&0&0
    \end{pmatrix}, \;
\rho_*(e_2) = \begin{pmatrix}
    0&1&0&0\\
    0&0&-1&0\\
    0&0&0&0\\
    0&0&0&0
    \end{pmatrix}, \;
\rho_*(e_3) = \begin{pmatrix}
    0&0&0&-1\\
    0&0&0&0\\
    0&0&0&0\\
    0&0&1&0
    \end{pmatrix}}.
\end{gather*}}

The general relations in \eqref{diagso13} are here
$\rho_*(e_j)\cdot \Lambda_\fm(u_i) - \Lambda_\fm(u_i) \cdot \rho_*(e_j)
- \Lambda_\fm([e_j,u_i]) = 0$, $i=1,\dotsc,4$, $j=1,2,3$.
In the cases $3.3^2(2,3,4)$ all matrices $\Lambda_\fm(u_i)$ vanish.
Then, the nonzero components of $R$ at $o$  are
$R_{23} = - \varepsilon \rho_*(e_2)$,
$R_{34}= - \varepsilon \rho_*(e_3)$.
In the case $3.3^2(4)$ all components of $R$ at $o$ vanish.
The 1st EYM eq.\ for the cases $3.3^2(2,3)$ is
{\small
\[
\begin{pmatrix}
    0&0&\bsl  a &0\\
    0&\bsl  a&0&0\\
    \bsl  a&0&\bsl b +  2\varepsilon  &0\\
    0&0&0& \bsl  a
\end{pmatrix}
=  \bsk
\begin{pmatrix}
    2/a&0&b/2a^2&0\\ \noalign{\smallskip}
    0&b/2a^2&0&0\\ \noalign{\smallskip}
    b/2a^2&0&b^2/2a^3&0\\ \noalign{\smallskip}
    0&0&0&b/2a^2
\end{pmatrix}.
\]}

Hence the cases $3.3^2(2,3)$ do not provide solutions of the EYM eqs.
Since one has $R=0$ for $3.3^2(4)$, this case neither furnish solutions.

\subsection{The cases ${\bf 3.5^1}(2)$, ${\bf 3.5^1}(3)$, ${\bf 3.5^1}(4)$.}

The common corresponding pair is
$\big( (\R^3 \times (\fso(1,2) ) \oplus  \R, \, \fso(1,2)\big)$.
The common nonzero brackets are
$[e_1,e_2] = 2e_2$, $[e_1,e_3] = -2e_3$, $[e_1,u_1]$ $= 2u_1$, $[e_1,u_3] = -2u_3$,
$[e_2,e_3] = e_1$, $[e_2,u_2] = u_1$, $[e_2,u_3] = -2u_2$, $[e_3,u_1] = 2u_2$,
$[e_3,u_2] = -u_3$, and the different ones are $[u_1,u_2] = \varepsilon e_2$,
$[u_1,u_3] = \varepsilon e_1$, $[u_2,u_3] = \varepsilon e_3$,
for $\varepsilon = 1, -1, 0$, respectively.
We have in the three cases $3.5^1(2,3,4)$ that
{\small
\begin{gather*}
g \equiv
\begin{pmatrix}
    0&0&2a&0\\
    0&a&0&0\\
    2a&0&0&0\\
    0&0&0&b
    \end{pmatrix}, \quad
\mathbf{r} \equiv
\begin{pmatrix}
0&0&-4\varepsilon&0 \\ \noalign{\smallskip}
0&-2\varepsilon &0&0 \\ \noalign{\smallskip}
-4\varepsilon&0&0&0 \\ \noalign{\smallskip}
0&0&0&0
\end{pmatrix},
\quad \bs = -\dfrac{6\varepsilon}{a},\\
{\setlength\arraycolsep{3pt}
\rho_*(e_1) =
\begin{pmatrix}
    2&0&0&0\\
    0&0&0&0\\
    0&0&-2&0\\
    0&0&0&0
    \end{pmatrix},\;
\rho_*(e_2) =
\begin{pmatrix}
    0&1&0&0\\
    0&0&-2&0\\
    0&0&0&0\\
    0&0&0&0
    \end{pmatrix}, \;
\rho_*(e_3) = \begin{pmatrix}
    0&0&0&0\\
    2&0&0&0\\
    0&-1&0&0\\
    0&0&0&0
    \end{pmatrix}}.
\end{gather*}}

All matrices $\Lambda_\fm(u_i)$, $i=1,\dotsc,4$, vanish.
Then, the nonnull components of $R$ at $o$  are
$R_{12} = -\varepsilon \rho_*(e_2)$,
$R_{13} = -\varepsilon \rho_*(e_1)$,
$R_{23} = -\varepsilon \rho_*(e_3)$ and the 1st EYM eq.\ is
{\small
\[
 {\setlength\arraycolsep{3pt}
\begin{pmatrix}
    0&0&2\bsl  a +2\varepsilon&0\\ \noalign{\medskip}
    0&\bsl  a +\varepsilon&0&0\\  \noalign{\medskip}
    2\bsl  a +2\varepsilon&0&0&0\\  \noalign{\medskip}
    0&0&0& \bsl b +(3\varepsilon b/a)
\end{pmatrix}
 = \bsk
\begin{pmatrix}
    1/a&0&-1/a&0\\ \noalign{\medskip}
    0&1/2a&0&0\\ \noalign{\medskip}
    -1/a&0&1/a&0\\  \noalign{\medskip}
    0&0&0&b/2a^2
\end{pmatrix}}.
\]}

Hence the cases $3.5^1(2,3)$ do not furnish any solution of the 1st EYM eq.
Neither the case $3.5^1(4)$.

\subsection{The cases ${\bf 3.5^2}(2)$, ${\bf 3.5^2}(3)$ ${\bf 3.5^2}(4)$.}
The corresponding pairs are $\big(\fso(1,3) \oplus  \R,\, \fso(3)\big)$,
$\big(\fso(4) \oplus  \R, \, \fso(3)\big)$ and
$\big((\R^3 \times \fso(3) ) \oplus  \R, \, \fso(3)\big)$,
respectively.

The common nonzero brackets are
$[e_1,e_2] = -e_3$, $[e_1,e_3] = e_2$, $[e_1,u_1] = -u_2$, $[e_1,u_2] = u_1$, $[e_2,e_3] = -e_1$,
$[e_2,u_1] = -u_3$, $[e_2,u_3] = u_1$, $[e_3,u_2] = -u_3$, $[e_3,u_3] = u_2$,
and the different ones are $[u_1,u_2] = \varepsilon e_1$, $[u_1,u_3] = \varepsilon e_2$, $[u_2,u_3] = \varepsilon e_3$,
for $\varepsilon = 1, -1, 0$, respectively.

We have in the three cases that
{\small
\begin{gather*}
g \equiv \diag(a,a,a,b), \; ab<0, \quad
\mathbf{r} \equiv \diag(-2\varepsilon, -2\varepsilon, -2\varepsilon,0), \quad
\bs = -6\varepsilon/a, \\
\rho_*(e_1) = \setlength\arraycolsep{3pt}
\begin{pmatrix}
    0&1&0&0\\
    -1&0&0&0\\
    0&0&0&0\\
    0&0&0&0
    \end{pmatrix}, \;
\rho_*(e_2) = \begin{pmatrix}
    0&0&1&0\\
    0&0&0&0\\
    -1&0&0&0\\
    0&0&0&0
    \end{pmatrix}, \;
\rho_*(e_3) = \begin{pmatrix}
    0&0&0&0\\
    0&0&1&0\\
    0&-1&0&0\\
    0&0&0&0
    \end{pmatrix}.
\end{gather*}}

All matrices $\Lambda_\fm(u_i)$, $i=1,\dotsc,4$, vanish.
Then, the nonnull components of $R$ at $o$ are
$R_{12} = - \varepsilon \rho_*(e_1)$,
$R_{13} = - \varepsilon \rho_*(e_2)$,
$R_{23} = -\varepsilon \rho_*(e_3)$ and the 1st EYM eq.\ is
{\small
\[ {\setlength\arraycolsep{3pt}
\begin{pmatrix}
\bsl a +\varepsilon &0&0&0 \\
0 &\bsl a +\varepsilon &0&0 \\
0&0&\bsl a +\varepsilon &0 \\
0&0&0& \bsl b + (3\varepsilon b/a)
\end{pmatrix} =
\bsk
\begin{pmatrix}
    1/2a&0&0&0\\ \noalign{\smallskip}
    0&1/2a&0&0\\ \noalign{\smallskip}
    0&0&1/2a&0\\ \noalign{\smallskip}
    0&0&0&-3b/2a^2
\end{pmatrix}},
\]}

\noindent with solutions
$\bsl=-3\varepsilon/2a$, $\bsk=-\varepsilon a$, $a\neq 0$.

In the case $3.5^2(4)$, one has $T \equiv 0$. Hence, according to Remark \ref{ntemt},
this case does not furnish any solution of the 1st EYM eq.

\subsection{The case ${\bf 4.1^2}(1)$.}
The corresponding pair is
{\small
\[
\left( \R^4 \times \left\{\begin{pmatrix}
x & z & 0 & -t \\
0 & 0 & -z & -y \\
0 & 0 &-x & 0 \\
0& y & t &0
 \end{pmatrix}\right\} ,\;
\left\{\begin{pmatrix}
x & z & 0 & -t \\
0 & 0 & -z & -y \\
0 & 0 &-x & 0 \\
0& y & t &0
 \end{pmatrix}\right\}\right),\; x,y,z,t \in \R.
\]}

The nonzero brackets are
$[e_1,e_3] = e_3$, $[e_1,e_4] = e_4$, $[e_1,u_1] = u_1$, $[e_1,u_3] = -  u_3$, $[e_2,e_3] = -e_4$, $[e_2,e_4] = e_3$, $[e_2,u_2] = u_4$, $[e_2,u_4] = -u_2$, $[e_3,u_2] = u_1$, $[e_3,u_3] =- u_2$, $[e_4,u_3] = u_4$, $[e_4,u_4]= -u_1$.

In this case, one has $T \equiv 0$. Hence, according to Remark \ref{ntemt},
this case does not furnish any solution of the 1st EYM eq.

\subsection{The cases ${\bf 6.1^3}(1)$, ${\bf 6.1^3}(2)$ and ${\bf 6.1^3}(3)$.}
The corresponding pairs are
$\big(\fso(3,2), \,\fso(1,3\big)$, $\big(\fso(1,4), \, \fso(1,3)\big)$ and
$\big( \R^4 \times \fso(1,3) , \, \fso(1,3)\big)$,
respectively.

The common nonzero brackets are
{\small
\begin{gather*}
[e_1,e_2] = -e_4, \;\; [e_1,e_3] = -e_5, \;\; [e_1,e_4] = e_2, \;\; [e_1,e_5] = e_3, \; [e_1,u_1] = -u_2, \\
[e_1,u_2] = u_1, \;\;
[e_2,e_3] = -e_6, \;\; [e_2,e_4] = -e_1, \;\; [e_2,e_6] = e_3, \;\; [e_2,u_1] =-u_3, \\
[e_2,u_3] = u_1, \quad
[e_3,e_5] = e_1,\quad [e_3,e_6] = e_2, \quad [e_3,u_1] = u_4, \quad [e_3,u_4] =u_1, \\
[e_4,e_5] = -e_6, \;\; [e_4,e_6] = e_5, \;\; [e_4,u_2] = -u_3, \;\; [e_4,u_3] = u_2, \quad
[e_5,e_6] = e_4, \\ [e_5,u_2] = u_4, \quad [e_5,u_4] = u_2, \quad [e_6,u_3] = u_4, \quad
[e_6,u_4] = u_3,
\end{gather*}}

\noindent
and the different ones are
$[u_1,u_2] = \varepsilon  e_1$, $[u_1,u_3] = \varepsilon  e_2$, $[u_1,u_4] = - \varepsilon e_3$,
$[u_2,u_3] = \varepsilon e_4$, $[u_2,u_4] = -\varepsilon e_5$, $[u_3,u_4] = -\varepsilon e_6$,
for $\varepsilon = 1, -1,0$, respectively.

One has in all cases $6.1^3(1,2,3)$ that
{\small
\begin{gather*}
g \equiv \diag(-a,-a,-a,a), \quad
\mathbf{r} \equiv \big(-3\varepsilon, -\varepsilon, -3\varepsilon, 3\varepsilon\big), \quad
\bs = 10\,\varepsilon/a,
\end{gather*}
\begin{gather*}
\rho_*(e_1) = {\setlength\arraycolsep{3pt}
\begin{pmatrix}
    0&1&0&0\\
    -1&0&0&0\\
    0&0&0&0\\
    0&0&0&0
    \end{pmatrix}, \;\;
\rho_*(e_2) =
\begin{pmatrix}
    0&0&1&0\\
    0&0&0&0\\
    -1&0&0&0\\
    0&0&0&0
    \end{pmatrix}}, \;\;
\rho_*(e_3) =
\begin{pmatrix}
    0&0&0&1\\
    0&0&0&0\\
    0&0&0&0\\
    1&0&0&0
    \end{pmatrix},\\
 {\setlength\arraycolsep{3pt}
\rho_*(e_4) = \begin{pmatrix}
    0&0&0&0\\
    0&0&1&0\\
    0&-1&0&0\\
    0&0&0&0
    \end{pmatrix}, \;
\rho_*(e_5) = \begin{pmatrix}
    0&0&0&0\\
    0&0&0&1\\
    0&0&0&0\\
    0&1&0&0
    \end{pmatrix}, \;
\rho_*(e_6) = \begin{pmatrix}
    0&0&0&0\\
    0&0&0&0\\
    0&0&0&1\\
    0&0&1&0
    \end{pmatrix}.}
\end{gather*}}

The general relations in \eqref{diagso13} are here
$\rho_*(e_j)\circ \Lambda_\fm(u_i) - \Lambda_\fm(u_i) \circ \rho_*(e_j)
- \Lambda_\fm([e_j,u_i]) = 0$, $i=1,\dotsc,4$, $j=1,\dotsc,6$;
thus in the cases $6.1^3(1,2,3)$ all matrices $\Lambda_\fm(u_i)$ vanish.
Then, the components of $R$ at $o$  are
$R_{12} = -\varepsilon\rho_*(e_1)$ ,
$R_{13} = -\varepsilon\rho_*(e_2)$,
$R_{14}= \varepsilon\rho_*(e_3)$,
$R_{23} = -\varepsilon\rho_*(e_4)$,
$R_{24} = \varepsilon\rho_*(e_5)$,
$R_{34}= \varepsilon\rho_*(e_6)$
and the 1st EYM equation reads
{\small
\[
\diag\big(-\bsl  a +2\varepsilon,- \bsl  a +4\varepsilon, -\bsl  a +2\varepsilon, \bsl  a -2\varepsilon\big)
= \bsk \, \varepsilon^2 \diag\Big(-\frac1a,\, -\frac1a,\,-\frac1a,\,-\frac3a\Big).
\]}

\noindent Hence, these cases do not provide any solution of the 1st EYM eq.

\subsection{The second EYM equation}

For those cases satisfying the first EYM equation, we now explore the behaviour of the second EYM equation.

\begin{prop}
Every connection $\Lambda$ satisfying the first Einstein-Yang-Mills equation in a Lorentzian four-codimensional symmetric pair $(\mathfrak{k},\mathfrak{h})$ also satisfies the second equation, that is, it is a Yang-Mills connection.
\end{prop}

\begin{proof}
We use the formula (see \cite[p.\ 381; and (3.5)]{EGH})
\begin{align}
\big(D(\Lstar R) \big)(X,Y,Z)  & = (\rd (\Lstar R))(X,Y,Z) \label{eym3} \\
     & \quad - (\Lstar R)(\Lambda_\fm(X)Y,Z)
                  - (\Lstar R)(Y, \Lambda_\fm(X)Z), \; \; X,Y,Z \in \fm. \notag
\end{align}
The first summand vanishes for symmetric spaces. In fact, let
$\{\theta^i\}$, $i=1,\dotsc,4$, be the basis dual to the orthonormal basis $\{u_j\}$, $j=1, \dotsc, 4$; that is, $\theta^i(u_j) = \delta^i_j$. Then the tensor $\rd(\Lstar R)$ can be expressed locally as a sum of $\fso(1,3)$-valued $2$-forms, each of its summands containing a factor of type $\rd \theta^i(u_j,u_k) = - \theta^i([u_j, u_k]_\fm)$. Now, one always has $[u_j,u_k]_\fm = 0$ for symmetric spaces. So \eqref{eym3} reduces to
\begin{equation}
\label{eym4}
\big(D (\Lstar R)\big)(u_i,u_j,u_k)  =
- (\Lstar R)(\Lambda_\fm(u_i)u_j,u_k) - (\Lstar R)(u_j, \Lambda_\fm(u_i) u_k).
\end{equation}
If all $\Lambda_\fm({u_j}) = 0$,  the Yang-Mills equation is satisfied. This happens for $3.5^2(2)$ and $3.5^2(3)$.
For the remaining eight cases, since the nonvanishing  $\Lambda_\fm(u_j)$ are those for $u_2,u_4$,
the right-hand side of \eqref{eym4} is possibly nonnull only in the two right-hand sides
 $ - (\Lstar R) (\Lambda_\fm({u_ 2})\cdot , \cdot) - (\Lstar R) (\cdot, \Lambda_\fm({u_ 2}) \cdot)$,  $-(\Lstar R)(\Lambda_\fm({u_ 4})\cdot , \cdot) - (\Lstar R) (\cdot, \Lambda_\fm({u_ 4}) \cdot) $.

Now, we always have that, up to a coefficient,
\[
\Lstar R =
\begin{cases}
\theta^2 \wedge \theta^4 \otimes R_{13}& \text{for}\; 1.1^1(7), 1.1^2(9,10), \\
\theta^2 \wedge \theta^4 \otimes R_{13} + \theta^1 \wedge \theta^3 \otimes R_{24}
& \text{for}\; 2.1^2(1,2,3), \\
\theta^1 \wedge \theta^3 \otimes R_{24} & \text{for} \; 2.1^2(4,5).
\end{cases}
\]
The only possibly nonnull RHS in \eqref{eym4} would be, up to a sign,
\[
v_{24}(\Lstar R)(u_4,u_2) + v_{24}(\Lstar R)(u_2,u_4),\quad
z_{24}(\Lstar R)(u_4,u_2) + z_{24} (\Lstar R)(u_2,u_4),
\]
which obviously vanish.
\end{proof}

 The surprising fact that a connection in a symmetric space is a Yang-Mills connection under some conditions (in our case, the condition derived from the first EYM equation) is a result that is compatible with previous results in the literate (for example, see \cite{taf}). In any case, as we already mentioned in the Preliminaries, we can be sure that the study of the EYM equations with non-vanishing currents in the second equation does not make sense in our context.

We have proved the following result.
\begin{prop}
\label{pr32}
\!The Lorentzian four-codimensional symmetric pairs $(\fk,\fh)$
furnishing solutions of the first Einstein-Yang-Mills
equation of type $(\Lambda,\!g^{\fl}_{\mathrm{diag}})$, are the ten ones listed in \emph{Table \ref{ccc}},
with their respective cosmological and gravitational constants $\bsl $ and $\bsk $.
\end{prop}
\begin{table}[htb]
\caption{}
\label{ccc}
\begin{tabular}{|l|c|c|c|c|c|} \hline
               &{\bf Case}  & $\bsl$  & $\bsk$  &{\bf Conds.} \\ \hline \hline
$\dim \fl =1$  & $1.1^2(7)$ & $\vphantom{(\frac12)^{3^4}}-1/(2a)$ & $a$ & $a\neq0$ \\
               & $1.1^2(9)$ & $1/(2a)$ & $a$ & $a\neq0$ \\
               & $1.1^2(10)$ & $-1/(2a)$ & $-a$ & $a\neq0$ \\
$\dim \fl =2$  & $2.1^2(1)$ & $\dfrac{a-b}{2ab}$ & $\dfrac{ab(a+b)}{a^2+b^2}$ & $a,b\neq0, \,a\neq -b$ \\
               & $2.1^2(2)$ & $\vphantom{(\frac12)^{{\frac34}^{5^{6^7}}}}-\dfrac{a+b}{2ab}$ & $\vphantom{(\frac12)^{{\frac34}^{5^{6^{7^8}}}}}-\dfrac{ab(a-b)}{a^2+b^2}$ & $a,b\neq0,\,a\neq b$  \\
               & $2.1^2(3)$   & $\vphantom{(\frac12)^{3^4}}-1/(2a)$ & $a$ & $a\neq0$ \\
               & $2.1^2(4)$  & $1/(2b)$ & $b$ & $b\neq0$  \\
               & $2.1^2(5)$  & $-1/(2b)$ & $-b$ & $b\neq0$\\
$\dim \fl =3$  & $3.5^2(2)$ & $-3/(2a)$ & $-a$ & $a\neq0$  \\
                     & $3.5^2(3)$ & $3/(2a)$ & $a$ & $a\neq0$ \\  \hline
\end{tabular}
\end{table}}

\section{Four-dimensional Lorentzian symmetric spaces furnishing solutions of Einstein-Yang-Mills equations of type $(\Lambda,g^\fl_\mathrm{diag})$}
\label{secfiv}

Denote by: $S^n$ the $n$-dimensional Riemannian sphere; $H^n$
the $n$-dimensional hyperbolic Riemannian space;
$\widetilde{E(2)} = \widetilde{(\R^2 \times \SO(2))/\SO(2)}$ the universal covering of the
Euclidean group of the plane $\R^2$; and $E(1,1) = (\R^2 \times \SO(1,1))/$ $\SO(1,1)$
the Euclidean group of the plane $\R^{1,1}$. Then we have
\begin{thm}
\label{th41}
The four-dimensional connected simply-connected
indecomposable Lorentz\-ian symmetric
spaces with connected isotropy group, furnishing solutions of the Einstein-Yang-Mills
equations of type $(\Lambda,g^{\fl}_{\mathrm{diag}})$,
are the ten spaces listed in \emph{Table~\ref{ddd}}.
\begin{table}[htb]
\caption{}
\label{ddd}
\begin{tabular}{|l|l||l|l|} \hline
{\bf Case}           & {\bf Space}                                               & {\bf Case} &  {\bf Space}  \\ \hline \hline
$1.1^1(7)$          & $\widetilde{H^2_1} \times \R^2$             & $2.1^2(3)$          & $\widetilde{E(2)} \times \vphantom{{\frac17}^{1^{1^1}}}\widetilde{H^2_1}$    \\
$1.1^2(9)$          & $S^2 \times \R^2$                                   & $2.1^2(4)$          & $S^2 \times E(1,1)$       \\
$1.1^2(10)$        & $H^2 \times \R^2$                                   & $2.1^2(5)$          & $H^2 \times E(1,1)$      \\
$2.1^2(1)$          & $S^2 \times \widetilde{H^2_1}$               & $3.5^2(2)$  & $H^3 \times \R$      \\
$2.1^2(2)$          & $H^2 \times \vphantom{{\frac17}^{1^{1^1}}}\widetilde{H^2_1}$               & $3.5^2(3)$ & $S^3 \times \R$   \\ \hline
\end{tabular}
\end{table}
\end{thm}
\begin{proof}
Let $(\fk,\fh,B)$ be a triple,
where $(\fk,\fh)$ is a reductive pair of Lie algebras and $B$ an $\Ad(H)$-invariant symmetric nondegenerate
bilinear form on the $\fh$-module $\fk/\fh$. By a result of Mostow \cite[p.\ 614]{Mos},
when $\mathrm{codim}_{\fk}\, \fh \leq 4$, there exists a unique (up to
equivalence)
homogeneous manifold $M = K/H$ such that $M$ is connected and simply-connected, the isotropy
group $H$ is connected, and the manifold admits a unique invariant
pseudo-Riemannian metric $g$ corresponding to the bilinear form~$B$.

The Komrakov pairs furnishing some nontrivial solution of the EYM eqs.\ are (Proposition  \ref{pr32})
{\small
\begin{align*}
& \big(\fsl(2) \oplus \R^2,\, \fso(1,1)\big), \hspace{22mm} \big((\R^2 \times \fso(2)) \oplus \fsl(2,\R),\, \fso(2) \oplus \fso(1,1) \big), \\
& \big(\fsu(2) \oplus \R^2,\, \fso(2)\big), \hspace{24mm} \big(\fso(3)\oplus(\R^2\times\fso(1,1)),\,\fso(2)\oplus\fso(1,1)\big), \\
& \big(\fsl(2,\R) \oplus \R^2,\, \fso(2)\big), \hspace{21mm} \big(\fso(2,\R)\oplus(\R^2\times \fso(1,1)),\,\fso(2)\oplus\fso(1,1)\big), \\
& \big(\fso(3) \oplus \fsl(2,\R),\, \fso(2) \oplus \fso(1,1) \big),\hspace{4mm} \big(\fso(1,3) \oplus  \R,\, \fso(3)\big), \\
& \big(\fsl(2,\R) \oplus \fsl(2,\R),\,\fso(2)\oplus \fso(1,1) \big),\hspace{1mm} \big(\fso(4) \oplus  \R, \, \fso(3)\big).
\end{align*}}

\noindent Table \ref{ddd} then follows because the corresponding connected simply-connected Lorentzian symmetric manifolds with connected isotropy group are
{\small
\begin{align*}
& \big(\widetilde{\SO_0(2,1)/\SO_0(1,1)}\big) \times \R^2,
\hspace{3mm}
\big((\R^2 \times \SO(2) )/\SO(2)\big) \times \big(\widetilde{\SO_0(2,1)/\SO_0(1,1)}\big), \\
& \big(\SO(3)/\SO(2)\big) \times \R^2, \hspace{12mm}
 \big(\SO(3)/\SO(2)\big) \times \big((\R^2 \times \SO_0(1,1) )/\SO_0(1,1)\big),\\ & \big(\SO_0(2,1)/\SO(2)) \times \R^2, \hspace{8mm}
\big(\SO_0(1,2)/\SO(2)\big) \times \big((\R^2 \times \SO_0(1,1) )/\SO_0(1,1)\big), \\
& \big(\SO(3)/\SO(2)\big) \times \big(\widetilde{\SO_0(2,1)/\SO_0(1,1)}\big), \hspace{6mm}
\big(\SO_0(1,3)/\SO(3)\big) \times \R, \\
& \big(\SO_0(1,2)/\SO(2)\big) \times \big(\widetilde{\SO_0(2,1)/\SO_0(1,1)}\big), \hspace{2mm}
\big(\SO(4)/\SO(3)\big) \times \R,
\end{align*}}
respectively.
\end{proof}

\begin{rem} \em
{\bf (a)} We have restricted ourselves, only for the sake of brevity, to a diagonal metric on each Lie algebra~$\fl$. The study of the general case, that is, considering not necessarily diagonal metrics on~$\fl$, remains an open problem.

{\bf (b)} The vanishing of the currents above is related with our choice of spacetimes being symmetric. As far as we know, the abundant bibliography on Einstein-Yang-Mills fields provide few explicit examples of nonnull currents in the context of homogeneous spacetimes. Many times the authors start either supposing that the current is null
ot that the group involved is $\SU(3)$ or other similar compact group.
In turn, each of the spacetimes defined by a four-dimensional symmetric Lie group admitting a biinvariant metric (see Levichev \cite{Lev1}), has null current. Thus, we think that the study in this sense of other Komrakov's pairs might enlighten the influence of the different sets of nonzero brackets in the expression of the currents.
This remains an open problem.
One first step would be for instance the study of the spaces with precisely one nonnull bracket $[u_i,u_j]_\fm$, that is, the eighteen spaces
$1.1^1(3,4,5,6)$, $1.1^2(3,4,5,6,7,8)$, $1.4^1(13,14,18,19,20$, $21,22,23)$.
\end{rem}

\subsection*{Acknowledgment}
We are indebted to Professor G. Calvaruso for the initial work we did on the topic
and his recent valuable suggestions, which have contributed to improve the manuscript a great deal.


\begin{thebibliography}{99}

\bibitem{BinSniFis} E.~Binz, J.~\'Sniatycki and H.\,R. Fischer,
\textit{ Geometry of Classical Fields},
North-Holland Mathematics Studies, vol.\ 154, North-Holland, Dordrecht, 1998.

\bibitem{Ble} D.~Bleecker,
\textit{ Gauge Theories and Variational Principles},
Addison-Wesley, Reading, MA, 1981.

\bibitem{CahWal} M.~Cahen and N.~Wallach,
\textit{ Lorentzian symmetric spaces},
Bull.\ Amer.\ Math.\ Soc.\ \textbf{ 76} (1970), no.\ 3, 585--591.

\bibitem{CV} G.~Calvaruso and J.~Van der Veken,
\textit{ Totally geodesic and parallel hypersurfaces of four-dimensional oscillator groups}. Results Math.\ \textbf{ 64} (2013), no.\ 1--2, 135--153.

\bibitem{DGO} R.~Dur\'an, P.\,M.~Gadea and J.\,A.~Oubi\~na,
\textit{ Reductive decompositions and Einstein-Yang-Mills equations associated to, the oscillator group},
J. Math.\ Phys.\ \textbf{ 40} (1999), no.\ 7, 3490--3498.

\bibitem{EGH} T.~Eguchi, P.\,B.~Gilkey, A.\,J.~Hanson,
\textit{ Gravitation, gauge theories and differential geometry},
Phys.\ Rep.\ \textbf{ 66} (1980), no.\ 6, 213--393.

\bibitem{KobNom} S.~Kobayashi and K.~Nomizu,
\textit{ Foundations of Differential Geometry, I, II},
Interscience Publishers, New York, 1963, 1969.

\bibitem{Kom1} B.~Komrakov, Jr.,
\textit{ Four-dimensional pseudo-Riemannian homogeneous spaces. Classification of complex pairs, I},
Preprint University of Oslo, no.\ \textbf{ 34}, December 1993.

\bibitem{Kom2} B.~Komrakov, Jr.,
\textit{ Four-dimensional pseudo-Riemannian homogeneous spaces. Classification of real pairs},
Preprint University of Oslo, no.\ \textbf{ 32}, June 1995.

\bibitem{Kom3} B.~Komrakov, Jr.,
\textit{ Einstein-Maxwell equation on four-dimensional homogeneous spaces},
Lobachevskii J. Math.\ \textbf{ 8} (2001), 33--165.

\bibitem{Lev1} A.\,V. Levichev,
\textit{ Some symmetric spaces of the general theory of relativity as solutions of the Einstein-Yang-Mills equations.} (Russian),
Group-theoretic methods in physics, vol.\ 1 (J$\bar{\mathrm u}$rmala, 1985), 145--150, ``Nauka,'' Moscow, 1986.

\bibitem{Lev2} A.\,V. Levichev,
\textit{ Pseudo-Hermitian realization of the Minkowski world through DLF theory},
Phys.\ Scr.\ \textbf{ 83} (2011), no.\ 1, 1--9.

\bibitem{Mos} G.\,D. Mostow,
\textit{ The extensibility of local Lie groups of transformations and groups on surfaces},
Ann.\ of Math.\ (2) \textbf{ 52} (1950), 606--636.

\bibitem{NW} G.\,R. Nappi and E.~Witten,
\textit{ Wess-Zumino-Witten model based on a nonsemisimple group},
Phys.\ Rev.\ Lett.\ \textbf{ 71} (1993), no.\ 23, 3751--3753.

\bibitem{Str} R.\,F. Streater,
\textit{ The representations of the oscillator group},
Commun.\ Math.\ Phys.\ \textbf{ 4} (1967), 217--236.

\bibitem{taf} J. Tafel, \textit{Some solutions of the Einstein-Yang-Mills equations}, in Geometrical
and topological methods in gauge theories (Proc. Canad. Math. Soc. Summer Res. Inst., McGill Univ., Montreal, Quebec, 1979),
pp. 134--136, Lecture Notes in Phys., vol. 129, Springer, Berlin-New York, 1980, p. 134.


\end{thebibliography}
\end{document}